\documentclass[12pt]{amsart}
\usepackage{amsmath,amsfonts,amssymb,amsthm}
\usepackage[dvips]{graphicx}
\usepackage{subfigure}
\usepackage{color}
\usepackage[T1]{fontenc}
\usepackage[utf8]{inputenc}

\theoremstyle{definition} 
\newtheorem{defs}{Definition}[section]

\theoremstyle{plain} 
\newtheorem{prop}[defs]{Proposition}
\newtheorem{lem}[defs]{Lemma}
 
\newtheorem*{NNthm}{Theorem} 

\newtheorem{cor}[defs]{Corollary}

\newtheorem*{ThmA}{Theorem A}
\newtheorem*{CorB}{Corollary B}

\theoremstyle{remark} 
\newtheorem{remark}[defs]{Remark}
\newtheorem*{NNrem}{Remark}
\newtheorem{Example}[defs]{Example}
\newtheorem{exemple}[defs]{Example}

\newtheorem*{Reader's guide}{Reader's guide}

\newcommand{\R}{\mathbb{R}}

\newcommand{\N}{\mathbb{N}}
\newcommand{\C}{\mathbb{C}}

\newcommand{\Z}{\mathbb{Z}}
\newcommand{\CC}{\mathcal{C}}
\newcommand{\HH}{\mathcal{H}}
\newcommand{\QQ}{\mathcal{Q}}
\newcommand{\HCC}{\widehat{\mathcal{C}}}
\newcommand{\HHCC}{\mathcal{H}_{\mathcal{C}}}
\newcommand{\HZ}{\widehat{Z}}

\title[Rauzy classes]
{Classification of Rauzy classes in the moduli space of Abelian and quadratic differentials}
\author{Corentin Boissy}
\address{Aix Marseille Université, CNRS, LATP, UMR~7353, 13453 Marseille France  }
\email{corentin.boissy@univ-amu.fr}

\subjclass[2000]{Primary: 37E05. Secondary: 37D40}
\keywords{Interval exchange maps, Linear involutions, Rauzy classes, Quadratic
 differentials, Moduli spaces}

\date{\today}

\begin{document}

\begin{abstract}
We study relations between Rauzy classes coming from an interval exchange map and the corresponding connected components of strata of the moduli space of Abelian differentials. This gives a criterion to decide whether two permutations are in the same Rauzy class or not, without actually computing them. We prove a similar result for Rauzy classes corresponding to quadratic differentials.
\end{abstract}

\maketitle

\section*{Introduction}

Rauzy induction was first introduced as a tool to study the dynamics of interval exchange transformations \cite{Rauzy}. These mappings appear naturally as first return maps on a transverse segment, of the directional flow on a translation surface. 
The Veech construction presents translation surfaces as suspensions over interval exchange maps, and extends the Rauzy induction to these suspensions \cite{Veech82}. This provides a powerful tool in the study of the Teichmüller geodesic flow and was widely studied in the last 30 years. 

An interval exchange map is encoded by a permutation and a continuous datum. A 
Rauzy class is a minimal subset of irreducible permutations which is invariant by the two combinatorial operations associated to the Rauzy induction. The Veech construction enables us to associate to a Rauzy class a connected component of the moduli space of Abelian differentials with prescribed singularities.
Such connected components are in one-to-one correspondence with the \emph{extended} Rauzy classes, which are unions of Rauzy classes and are defined by adding a third combinatorial operation. 
Historically, these extended Rauzy classes were used to prove the nonconnectedness of some strata in low genera \cite{Veech90}, before Kontsevich and Zorich performed the complete classification \cite{Kontsevich:Zorich}.
\medskip

One can also consider first return maps of the vertical foliation on transverse segments for flat surfaces defined by a quadratic differential on a Riemann surface.  We obtain a particular case of \emph{linear involutions}, that were defined by Danthony and Nogueira \cite{DaNo} as first return maps of measured foliations on surfaces. In this paper, we speak only of linear involutions corresponding to quadratic differentials. As before, a linear involution is encoded by a combinatorial datum, the \emph{generalized permutation} and a continuous datum. 
For linear involutions with irreducible generalized permutations, we can generalize the Veech construction and Rauzy classes \cite{Boissy:Lanneau}.

In this paper, we give a precise relation between Rauzy classes and the connected components of the moduli space of Abelian or quadratic differentials. We prove the following:

\begin{ThmA}
Let $\QQ$ be a stratum in the moduli space of Abelian differentials or in the moduli space of quadratic differentials. Let $r$ be the number of distinct orders of singularities of an element of~$\QQ$. For any connected component $\CC$ of $\QQ$, there are exactly $r$ distinct Rauzy classes that correspond to this connected component.
\end{ThmA}

This gives a positive answer to Conjecture~2 stated in \cite{Zorich:jmd}.
Note that in the previous theorem, $r$ is not the number of singularities: for instance, in the stratum that consists of translation surfaces with two singularities of degree $1$ (\emph{i.e.} the stratum $\HH(1,1)$), we have $r=1$.
\medskip

Theorem~A will be obtained as a direct combination of Propositions~\ref{RC:Cm} and~\ref{CC:Cm:ab} for the case of Abelian differentials, and Propositions~\ref{RC:Cm} and \ref{CC:Cm:q} for the case of quadratic differentials.
\medskip

A flat surface obtained from a permutation or a generalized permutation $\pi$ using the Veech construction admits a marked singularity. The order of this singularity $\alpha(\pi)$ is preserved by the Rauzy induction, and we can therefore associate to a Rauzy class an integer, which is the order of a singularity in the corresponding stratum. Hence, a corollary of Theorem~A is the following criteria:

\begin{CorB}
Let $\pi_1$ and $\pi_2$ be two irreducible permutations or generalized permutations. They are in the same Rauzy class if and only if they correspond to the same connected component and $\alpha(\pi_1)=\alpha(\pi_2)$.
\end{CorB}

See Appendix A for further comments concerning this corollary.

\subsubsection*{Reader's guide}
In section 1, we recall the definition and some facts about flat surfaces. In particular, we present the ``breaking up singularities" surgeries on flat surfaces that will be an essential tool for the proof of the main result. Note that the surgery presented in section~\ref{sec:non:loc} is more technical and can be skipped in a first reading.\\
In section 2, we recall the definitions about interval exchange, linear involutions, and Rauzy classes.\\
In section 3, we show that there is a one-to-one correspondence between  Rauzy classes and connected components of the moduli space of flat surfaces with a marked singularity. This is Proposition~\ref{RC:Cm} .
\\
In section 4, we classify the connected components of the moduli space of flat surfaces with a marked singularity. This will correspond to Proposition~\ref{CC:Cm:ab} for Abelian differential and Proposition~\ref{CC:Cm:q} for quadratic differentials. Then, Theorem~A will follow directly from the main results of Section~3 and Section~4.

\subsubsection*{Acknowledgments}
 I  thank Anton Zorich, Pascal Hubert and Erwan Lanneau for encouraging me to write this paper, and for many discussions. I am gratefull to the Max-Planck-Institut at Bonn for its hospitality. I also thank the  anonymous referee for comments and remarks.

\section{Flat surfaces} 
\subsection{Definition}

A {\it flat surface} is a real, compact, connected surface of genus $g$ equipped with a flat metric with isolated conical singularities and such that the linear holonomy group belongs to $\Z / 2\Z$. Here holonomy means that the parallel transport of a vector along any loop brings the vector back to itself or to its opposite. This implies that all cone angles are integer multiples of $\pi$. We also fix a choice of a parallel line field in the complement of the conical singularities. This parallel line field will be usually referred as \emph{the vertical direction}.
 Equivalently a flat surface is a triple $(S,\mathcal U,\Sigma)$ such that $S$ is a topological compact connected surface, $\Sigma$ is a finite subset of $S$ (whose elements are called {\em singularities}) and $\mathcal U = \{(U_i,z_i)\}$ is an atlas of $S \setminus \Sigma$ such that the transition maps $z_j \circ z_i^{-1} : z_i(U_i\cap U_j) \rightarrow z_j(U_i\cap U_j)$ are translations or half-turns: $z_i = \pm z_j + c$, and for each $s\in \Sigma$, there is a neighborhood of $s$ isometric to a Euclidean cone. Therefore, we get a {\it quadratic differential} defined locally in the coordinates $z_i$ by
the formula $q=d z_i^2$. This quadratic differential extends to the points of $\Sigma$ to zeroes, simple poles or marked points (see~\cite{Masur:Tabachnikov}). Slightly abusing  vocabulary, a pole will be referred to as a zero of order $-1$, and a marked point will be referred to as a zero of order $0$. Then, a zero of \emph{order} $k\geq -1$ corresponds to a conical singularity of angle $(k+2)\pi$.  

Observe that the linear holonomy given by the flat metric is trivial if and only if there exists a sub-atlas such that all transition functions are translations or equivalently if the quadratic differential $q$ is the global square of an Abelian differential. We will then say that $S$ is a translation surface. In this case, we can choose a parallel vector field instead of a parallel line field, which is equivalent in fixing a square root $\omega$  of $q$. Also, a zero of \emph{degree} $k\geq 0$ of $\omega$ corresponds to a conical singularity of angle $(k+1)2\pi$. 

When a flat surface is not a translation surface, \emph{i.e.} if the corresponding quadratic differential is not the square of an Abelian differential, we oftently use the terminology \emph{half-translation surfaces}, since the change of coordinates are either translations or half-turns.

Following a convention of Masur and Zorich (see \cite{Masur:Zorich}, section~5.2), we will speak of the \emph{degree} of a singularity in a translation surface, and of the \emph{order} of a singularity in half-translation surface, since one of them refer to a zero of an Abelian differential and the other to a quadratic differential.

\begin{exemple}
 Consider a polygon whose sides come by pairs, and such that, for each pair, the corresponding sides are parallel and have the same length. We identify each pair of sides by a translation or a half-turn so that it preserves the orientation of the polygon. We obtain a flat surface, which is a translation surface if and only if all the identifications are done by translation. One can show that any flat surface can be represented by such a polygon (see \cite{B2}, Section~2).
\end{exemple}

A saddle connection is a geodesic segment (or geodesic loop) joining two singularities (or a singularity to itself) with no singularities in its interior. Even if $q$ is not globally a square of an Abelian differential, we can find a square root of $q$ along the interior of any saddle connection. Integrating $q$ along the saddle connection we get a complex number (defined up to multiplication by $-1$). Considered as a planar vector, this complex number represents the affine holonomy vector along the saddle connection. In particular, its Euclidean length is the modulus of its holonomy vector. 
\smallskip

\subsection{Moduli spaces}

For $g  \geq 0$, we define  the moduli space  of quadratic differentials
$\mathcal{Q}_g$ as the moduli space of pairs $(X,q)$ where $X$ is a genus
$g$  (compact,  connected)  Riemann  surface  and  $q$ a  non-zero quadratic differential $X$.  The term
moduli  space  means that  we  identify  the  points $(X,q)$  and
$(X',q')$   if  there   exists  an   analytic   isomorphism  $f:X
\rightarrow X'$ such that $f^\ast q'=q$.  Equivalently, in terms of polygon representations, two flat surfaces are identified in the moduli space of quadratic differentials if and only if the corresponding polygons can be obtained from each other by some finite number of ``cutting and gluing'',  preserving the identifications. The moduli space of Abelian differentials $ \mathcal{H}_{g}$, for $g\geq 1$ is defined in a analogous way.

We can associate to a quadratic differential the set with multiplicities $\{k_1^{\alpha_1},\ldots,k_r^{\alpha_r}\}$ of orders $\{k_1,\ldots,k_r\}$ of its poles and zeros, where $k_i\neq k_j$ for $i\neq j$, and  $k_i\geq -1$ and $\alpha_i\geq 1$ is the multiplicity of $k_i$.  
 The Gauss--Bonnet formula asserts that $\sum_i \alpha_i k_i=4g-4$. Conversely, if we fix a set with multiplicities $\{k_1^{\alpha_1},\dots,k_r^{\alpha_r}\}$ of integers, greater than or equal to $-1$ satisfying the previous equality, we denote by $\mathcal{Q}(k_1^{\alpha_1},\ldots,k_r^{\alpha_r})$ the  moduli space of quadratic differentials which are not globally squares of Abelian differentials, and which have $\{k_1^{\alpha_1},\ldots,k_r^{\alpha_r}\}$ as orders of poles and zeros. By a result of Masur and Smilie \cite{Masur:Smillie},  this space is nonempty except for $\mathcal{Q}(\emptyset)$, $ \mathcal{Q}(3,1)$, $ \mathcal{Q}(4)$ and $ \mathcal{Q}(-1,1)$. In the nonempty case,
 it is well known that  $\mathcal{Q}(k_1^{\alpha_1},\ldots,k_r^{\alpha_r})$ is a complex analytic orbifold, which is usually called a \emph{stratum} of the moduli space of quadratic differentials on a Riemann surface of genus $g$. 
  In a similar way, we denote by $\mathcal{H}(n_1^{\alpha_1},\dots,n_r^{\alpha_r})$ the moduli space of Abelian differentials 
having zeroes of degree $\{n_1^{\alpha_1},\ldots,n_r^{\alpha_r}\}$, where $n_i\geq 0$ and $\sum_{i=1}^r \alpha_i n_i=2g-2$.

There is a natural action of $\textrm{SL}_2(\mathbb{R})$ on each strata: let $(U_i,\phi_i)_{i\in I}$ be an atlas of flat coordinates of $S$, with~$U_i$ open subset of $S$ and $\phi_i(U_i)\subset \mathbb{R}^2$. An atlas of $A.S$ is given by $(U_i,A\circ \phi_i)_{i\in I}$. The action of the diagonal subgroup of $\textrm{SL}_2(\mathbb{R})$ is called the Teichmüller geodesic flow. In order to specify notations, we denote by $g_t$ the matrix 
$\left(\begin{smallmatrix} e^{t/2}& 0 \\ 0& e^{-t/2} \end{smallmatrix}\right)$.
\smallskip

Local coordinates for a stratum of Abelian differentials are obtained by integrating the holomorphic 1--form along a basis of the relative homology $H_1(S,\Sigma;\mathbb{Z})$, where $\Sigma$ denotes the set of conical singularities of $S$. Equivalently, this means that local coordinates are defined by the relative cohomology $H^1(S,\Sigma;\mathbb{C})$.

Local coordinates in a stratum of quadratic differentials are obtained in the following way (see for instance \cite{Douady:Hubbard}): one can naturally associate to a quadratic differential $(S,q)\in \mathcal{Q}(k_1^{\alpha_1},\ldots,k_r^{\alpha_r})$ a double cover $p: \widehat{S}\rightarrow S$ such that $p^* q$ is the square of an Abelian differential $\omega$. Let $\widehat{\Sigma}=p^{-1}(\Sigma)$. The surface $\widehat{S}$ admits a natural involution $\tau$, that induces on the relative homology $H_1(\widehat{S},\widehat{\Sigma};\mathbb{Z})$ an involution $\tau^*$. It decomposes $H_1(\widehat{S}, \widehat{\Sigma};\mathbb{Z})$ into an invariant subspace $H_1^+(\widehat{S},\widehat{\Sigma};\mathbb{Z})$ and an anti-invariant subspace $H_1^-(\widehat{S},\widehat{\Sigma};\mathbb{Z})$. 
Then local coordinates for a stratum of quadratic differential are obtained by integrating $\omega$ along a basis of  $H_1^-(\widehat{S},\widehat{\Sigma};\mathbb{Z})$.

\subsection{Connected components of the moduli space of Abelian differentials}
Here, we recall the classification of the connected components of the strata of the moduli space of Abelian differentials, due to Kontsevich and Zorich \cite{Kontsevich:Zorich}.

\begin{defs}
A flat surface $S$ is called hyperelliptic if there exists an orientation preserving involution $\tau$ which preserves the flat metric such that $S/\tau$ is a (flat) sphere.
\end{defs}

Sometimes, a connected component of a stratum consists only of hyperelliptic flat surfaces. In this situation it is called a hyperelliptic connected component. 
\medskip

Let $\gamma$ be a smooth curve in $S$ that does not contains any singularity. We parametrize $\gamma$ by arc length. In a translation surface, there is a natural identification between $\C$ and the tangent space of a regular point.  Hence,  one can identify $\gamma'$ to a closed path in the unit circle of~$\mathbb{C}$, \emph{e.g.} using the Gauss map, and compute its index that we denote by $Ind(\gamma)$. 

\begin{defs}[Kontsevich-Zorich]
Let $(\alpha_i,\beta_i)_{i\in \{1,\ldots,g\}}$ be a collection of paths representing a symplectic basis for the homology $H_1(S;\mathbb{Z})$.  We define the \emph{parity of the spin structure} of $S$ to be:
$$\sum_{i=1}^g \left(Ind(\alpha_i)+1\right)\left(Ind(\beta_i)+1\right) \mod 2.$$
\end{defs}

If all the singularities of the surface are of even degree, one can show that the parity of the spin structure does not depend on the choice of the paths and is an invariant of the connected component of the corresponding stratum. Now we can state the classification of these connected components.

\begin{NNthm}[Kontsevich-Zorich]
Let $\HH=\HH(k_1^{\alpha_1},\ldots,k_r^{\alpha_r})$ be a stratum in the moduli space of Abelian differentials, with $k_i\neq k_j$ for $i\neq j$, and with $k_i> 0$ and $\alpha_i>0$ for all $i$. Let $g$ be the corresponding genus. The stratum $\HH$ admits one, two, or three connected components according to the following rules:
\begin{enumerate}
\item If $\HH=\HH(2g-2)$ or $\HH(g-1,g-1)$, then $\HH$ contains one hyperelliptic connected component. If $g=2$, this component is the whole stratum, and if $g=3$, there is exactly one other connected component.
\item If $g\geq 4$ and if $k_1,\ldots,k_r$ are even, then there are exactly two connected components of $\HH$, with different parity of spin structures, and 
that are not hyperelliptic components.
\item In any other case, the stratum $\HH$ is connected.
\end{enumerate}
\end{NNthm}
Note that in the previous statement, the cases 1 and 2 can occur simultaneously. For instance, the stratum $\HH(6)$ has three connected components: one hyperelliptic, and two others that are distinguished by the parities of the corresponding spin structures.

\begin{remark}
The theorem  above is given for strata with no marked points. The classification for strata with marked points, \emph{i.e.} where we authorize $k_i=0$, is deduced in an obvious way.  
\end{remark}

\subsection{Connected components of the moduli space of quadratic differentials}
In this section, we recall the  classification of connected components of the strata in the moduli space of quadratic differentials, that will be needed (see \cite{Lanneau:hc,Lanneau:cc}).

\begin{NNthm}[E. Lanneau]
The hyperelliptic connected components are given by the following list:
\begin{enumerate}
\item The subset of surfaces in $\mathcal{Q}(k_1,k_1,k_2,k_2)$, that are a double covering of a surface in $\mathcal{Q}(k_1,k_2,-1^{s})$ ramified over $s$ poles. Here $k_1$ and $k_2$ are odd, $k_1\geq -1$ and $k_2\geq 1$, and $k_1+k_2- s=-4$.
\item The subset of surfaces in $\mathcal{Q}(k_1,k_1,2k_2+2)$, that are a double covering of a surface in $\mathcal{Q}(k_1,k_2,-1^{s})$ ramified over $s$ poles and over the singularity of order $k_2$.
Here $k_1$ is odd and $k_2$ is even, $k_1\geq -1$ and $k_2\geq 0$, and
$k_1+k_2-s=-4$.
\item The subset of surfaces in $\mathcal{Q}(2k_1+2,2k_2+2)$, that are a double covering of a surface in $\mathcal{Q}(k_1,k_2,-1^{s})$ ramified over all the singularities. 
Here $k_1$ and $k_2$ are even, $k_1\geq 0$ and $k_2\geq 0$, and $k_1+k_2-
s=-4$.
\end{enumerate}
\end{NNthm}

\begin{NNthm}[E. Lanneau]
In the moduli space of quadratic differentials, the nonconnected strata have two connected components and are in the following list (up to marked points):
\begin{itemize}
\item The strata that contain a hyperelliptic connected component, except the following ones, that are connected: $\QQ(-1,-1,-1,-1)$, $\QQ(-1,-1,1,1)$, $\QQ(-1,-1,2)$, $\QQ(1,1,1,1)$, $\QQ(1,1,2)$ and $\QQ(2,2)$.
\item The \emph{exceptionnal} strata $\QQ(-1,9)$, $\QQ(-1,3,6)$, and $\QQ(-1,3,3,3)$  and $\QQ(12)$.
\end{itemize}
\end{NNthm}

\subsection{Breaking up a singularity: local construction}\label{bzero:loc}
Here we describe a surgery, introduced by Eskin, Masur and Zorich (see \cite{EMZ}, Section~8.1) for Abelian differentials, that ``break up'' a singularity of degree $k_1+k_2\geq 2$ into two singularities of degree $k_1\geq 1$ and $k_2\geq 1$ respectively. This surgery is \emph{local}, since the metric is modified only in a neighborhood of the singularity of degree $k_1+k_2$.  The case $k_1=0$ or $k_2=0$ is trivial.  

\begin{figure}[htb]
\begin{center}
\includegraphics[width=360pt]{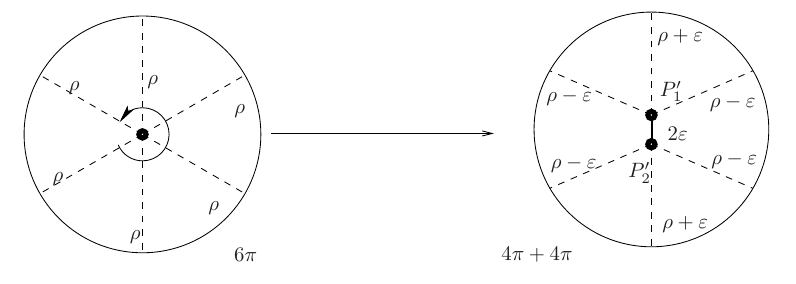}
\caption{Breaking up a zero, after Eskin, Masur and Zorich}
\label{bzero}
\end{center}
\end{figure}

We start from a singularity of degree $k_1+k_2$. A neighborhood of such singularity is obtained by gluing $(2k_1+2k_2+2)$ Euclidean half disks in a cyclic order. The singularity breaking procedure consists in changing continuously the way these half disks are glued together, as in Figure~\ref{bzero}. This breaks the singularity of degree $k_1+k_2$ into  singularities of degree $k_1$ and $k_2$ respectively, and with a small saddle connection joining them. 

Note that since the previous procedure purely local, it is also valid for quadratic differentials, as soon as we break up a singularity of even order into two singularities of even order. One can also in a similar way break up a singularity of odd order into a pair of singularities (see \cite{Masur:Zorich} for instance) although we will not need that case.  One can show that it is not possible to break a singularity of even order into two singularities of odd order by a local surgery. We need for this a nonlocal construction.

\subsection{Breaking up a singularity: nonlocal constructions} \label{sec:non:loc}
Here we describe a surgery, introduced by Masur and Zorich (see \cite{Masur:Zorich}, Section~6) for quadratic differentials, that ``break up'' a singularity of order $k_1+k_2$ into two singularities of order $k_1$ and $k_2$ respectively. It is valid for any $k_1,k_2\geq -1$, with $(k_1,k_2)\neq (-1,-1)$.  

\begin{figure}[htb]
\begin{center}
\includegraphics[width=250pt]{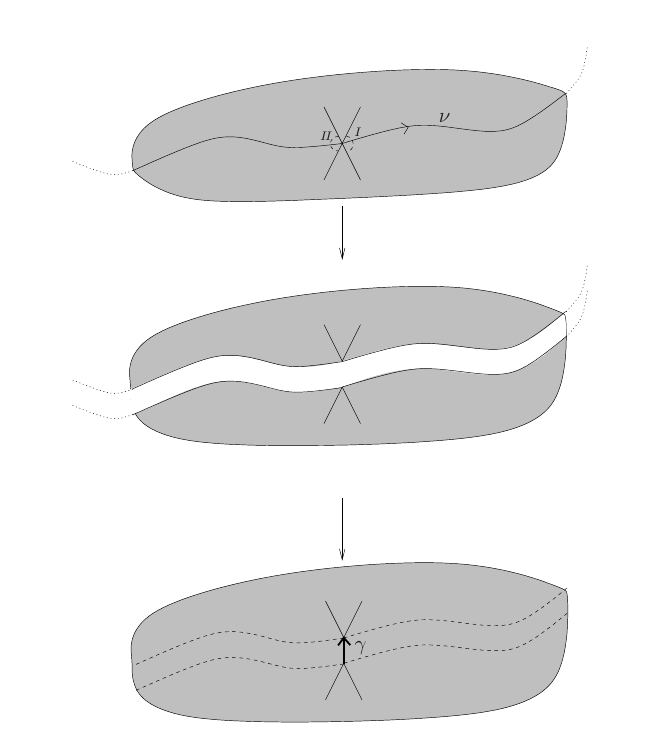}
\caption{Breaking up a zero of order two into two zeros of order one.}
\label{bzero:non:loc}
\end{center}
\end{figure}

We start from a surface $S_0$ with a singularity of order $k_1+k_2$, and other singularities of order $n_1,\ldots,n_s$. 
Consider an angular sector of angle $\pi$ between two consecutive vertical separatrices of $P$. We denote by $I$ this sector and by $I\!I$ the image of $I$ by a rotation of angle $(k_1+1)\pi$, and of center $P$.
Then, choose a closed path $\nu$ transverse to the vertical foliation that starts from the singularity $P$, sector~$I$ and ends at $P$, sector~$I\!I$. We also ask that the path $\nu$ does not intersect any singularity except $P$ in its end points. Then, we cut the surface along this path and paste in a ``curvilinear annulus'' with two opposite sides isometric to $\nu$, and with vertical height of length $\varepsilon$ (see Figure~\ref{bzero:non:loc}). We get a surface with singularities of order $k_1,k_2,n_1,\ldots,n_s$, with the same holonomy as $S_0$, and with a simple saddle connection $\gamma$ joining the two newly created singularities of order $k_1$ and $k_2$ . We denote this flat surface by $S=\Psi(S_0,\nu, \varepsilon)$. 
Similarly, we can perform the same construction, using the foliation $\mathcal{F}_\theta$ of angle $\theta$, and a path $\nu$ transverse to the foliation $\mathcal{F}_\theta$. We get a surface $\Psi_\theta(S_0,\nu,\varepsilon)$.

Note that giving an orientation to $\nu$ gives an orientation to $\gamma$ in the following way: $\nu$ defines a element $[\nu]$ in the homotopy group of $S\backslash\Sigma$, where $\Sigma$ is the set of conical singularities of $S$. The intersection number between $\gamma$ and $[\nu]$ is $\pm 1$ depending on the orientation of~$\gamma$. We then fix the orientation of $\gamma$ such that this intersection number is one. Then, we can consider $S=\Psi(S_0,\nu, \varepsilon)$ as an element of the moduli space of quadratic differentials with a \emph{marked} singularity by saying that the marked point of $S$ is the starting point of $\gamma$.

This construction was generalized by the author to polygonal curves in \cite{B2}, section~3. Such curve must still be transverse to the vertical foliation in a neighborhood of the singularity $P$ and must have nontrivial linear holonomy (if $k$ is odd). If $\nu$ is such path, then for $\varepsilon$ small enough, we get a surface $S=\Psi(S_0,\nu, \varepsilon)$ as described in the previous paragraph (by a surgery performed in a neighborhood of $\nu$). 
This new construction is more flexible and we have the following facts.

\begin{enumerate}
\item $\Psi(S_0,\nu, \varepsilon)$ depends continuously on $\varepsilon$ and on $S_0$.
\item If $\gamma\subset S$ is a vertical saddle connection joining two different singularities and is very small compared to any other saddle connection of $S$, then there exists a flat surface $S_0$ and $\nu_0\subset S_0$ such that $S=\Psi(S_0,\nu_0,\varepsilon)$ (see \cite{B2}, proof of Proposition~4.6).
\item The flat surface $\Psi(S_0,\nu_0,\varepsilon)$ does not change under small perturbations of $\nu_0$ (see \cite{B2}, Corollary~3.5).
\item Let $\nu_1$ be another path on $S_0$ that does not intersect any singularities except $P$ and starts and ends on sectors $I,I\!I$ of $P$ respectively. There exists $S_1$ in a neighborhood of $S_0$ such that $\Psi(S_0,\nu_1,\varepsilon)=\Psi(S_1,\nu_0,\varepsilon)$, and $S_1$ can be chosen arbitrarily close to $S_0$ as soon as $\varepsilon$ is small enough (\cite{B2}, proof of Lemma~4.5).
\end{enumerate}

\section{Rauzy classes}
\subsection{Interval exchange maps and linear involutions}
The first return map of the vertical flow of a translation surface on a horizontal open segment $X$ defines an \emph{interval exchange map}. That is, a one-to-one map from $X\backslash \{x_1,\dots,x_{d-1}\}$ to $X\backslash \{x_1',\ldots,x_{d-1}'\}$ which is an isometry and preserves the natural orientation of $X$. The relation between translation surfaces and interval exchange transformations has been widely studied in the last 25 years (see~\cite{Keane, Katok, Veech82,Masur82,Marmi:Moussa:Yoccoz, Avila:Gouezel:Yoccoz,Avila:Viana} etc\ldots).

\begin{figure}[htbp]
 \begin{center}
 \includegraphics{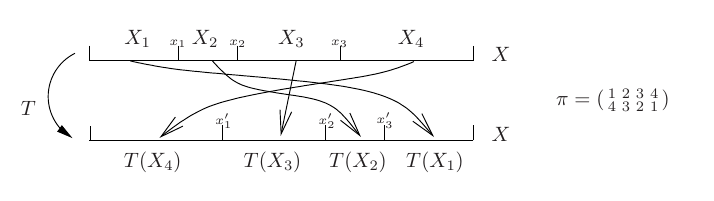}
 \caption{An interval exchange map and its corresponding permutation.}
 \end{center}
\end{figure}

We encode an interval exchange map $T$ in the following way: the set $X\backslash \{x_1,\ldots,x_{d-1}\}$ is a union of $d$ intervals that we label by $\{1,\ldots,d\}$ from the left to the right. The length of these intervals is then a vector $\lambda$ with positive entries. 
Applying the map $T$, the interval number $i$ is mapped to the interval number $\pi(i)$. This defines a permutation $\pi$ of $\{1,\ldots,d\}$. The vector $\lambda$ is called the continuous datum of $T$ and $\pi$ is called the combinatorial datum. We usually represent $\pi$ by a table of two lines:
\[
\pi=\left(\begin{array}{cccc}
1&2&\ldots&d \\ \pi^{-1}(1)&\pi^{-1}(2)&\ldots&\pi^{-1}(d)
\end{array}\right).
\]

 \medskip

The vertical foliation of a translation surface is a oriented measured foliation on a smooth oriented surface. A generalization of interval exchange maps for any measured foliation on a surface (oriented or not) was introduced by Danthony and Nogueira \cite{DaNo} as \emph{linear involution}. 
 The linear involutions corresponding to oriented flat surfaces with $\Z/2\Z$ linear holonomy were studied in detail by Lanneau and the author in \cite{Boissy:Lanneau}. 
 
 Let $X\subset S$ be an open horizontal segment. We choose on $X$ an orientation. This is equivalent to fix a ``left end'' on $X$, or to fix a ``positive vertical direction'' in a neighborhood of $X$. A linear involution must encode the successive intersections of $X$ with a vertical geodesic. It is done in the following way: we say that we are in $X\times \{0\}$ if the geodesic intersects $X$ in the positive direction and in $X\times \{1\}$ in the complementary case. Then, the first return map with this additional directional information gives a map from $X\times\{0,1\}$ to itself.

\begin{figure}[htbp]

\begin{center}



\includegraphics[width=8cm]{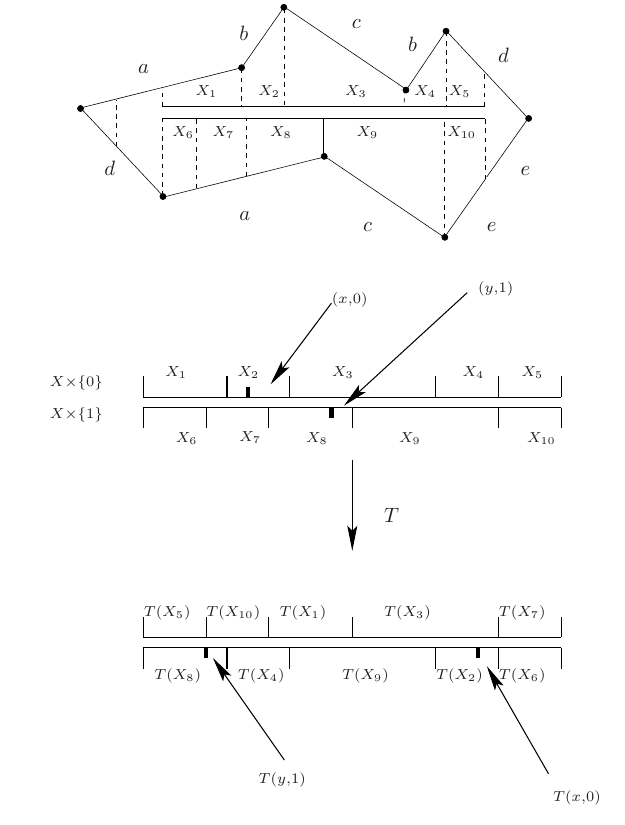}

 \caption{A linear involution associated to a measured foliation on a 
flat surface.}
 \label{ex:giem}
 \end{center}
\end{figure}

 \begin{defs}
 Let $f$ be the involution of $X\times\{0,1\}$ given by $f(x, \varepsilon)=(x,1-\varepsilon)$. A \emph{linear involution} is a map $T$, from $X\times\{0,1\}$ into itself, of the form $f\circ \tilde T$, where $\tilde T$ is an involution of $X\times\{0,1\}$ without fixed point, continuous except on a finite set of points $\Sigma_T$, and which preserves the Lebesgue measure. In this paper we will only consider linear involutions with the following additional condition: the derivative of $\tilde T$ is $-1$ at $(x,\varepsilon)$ if
$(x,\varepsilon)$ and $\tilde T(x,\varepsilon)$ belong to the same connected component, and $+1$ otherwise. 
\end{defs}
\medskip

On a flat surface, the first return map of the vertical foliation on a horizontal segment defines a linear involution. The fact that the underlying flat surface is oriented corresponds precisely to our additional condition. A linear involution such that $T(X\times \{0\})=X\times \{0\}$ (up to a finite subset) corresponds to an interval exchange map $T_0$, by restricting $T$ on $X\times \{0\}$ (note that the restriction of $T$ on $X\times \{1\}$ is naturally identified with $T_0^{-1}$). Therefore, we can identify the set of interval exchange maps with a subset of the linear involutions.

A linear involution is encoded by a combinatorial datum called \emph{generalized permutation} and by continuous data. 
This is done in the following way: 
 $X\times \{0\}\backslash \Sigma_T$ is a union of $l$ open intervals $X_1\sqcup \ldots \sqcup X_l$, where we assume by convention that $X_i$ is the interval at the place $i$, when counted from the left to the right. Similarly, $X\times \{1\}\backslash \Sigma_T$ is a union of $m$ open intervals $X_{l+1}\sqcup \ldots \sqcup X_{l+m}$. For all $i$, the image of $X_i$ by the map $\tilde{T}$ is a interval $X_j$, with $i\neq j$, hence $\tilde T$ induces an involution without fixed points on the set $\{1,\ldots, l+m\}$. To encode this involution, 
we attribute to each interval $X_i$ a symbol such that $X_i$ and $\tilde T(X_i)$ share the same symbol. Choosing the set of symbol to be $\{1,\ldots,d\}$, 
we get a two-to-one map $\pi:\{1,\ldots,l+m\}\rightarrow \{1,\ldots,d\}$, with $d=\frac{l+m}{2}$. Note that $\pi$ is not uniquely defined by $T$ since we can compose it on the left by any permutation of $\{1,\ldots,d\}$.

\begin{defs}
A \emph{generalized permutation} of type $(l,m)$, with $l+m=2d$, is a two-to-one map $\pi:\{1,\ldots,2d\}\rightarrow \{1,\ldots,d\}$. It is called \emph{reduced} if for each $k$, the first occurrence in $\{1,\ldots,l+m\}$ of the label $k\in\{1,\ldots,d\}$ is before the first occurrence of any label $k'>k$. 

We will usually represent such generalized permutation by a table of two lines of symbols, with each symbol appearing exactly two times.
$$
\pi=\left(\begin{array}{ccc}
\pi(1) & \dots& \pi(l) \\
\pi(l+1) & \dots & \pi(l+m)
\end{array}\right). 
$$
\end{defs}

In the table representation of a generalized permutation, a symbol might appear two times in a line, and zero time in the other line. Therefore, we do not necessarily have $l=m$.  A linear involution defines a reduced generalized permutation by the previous construction in a unique way. 

\begin{exemple}
The reduced generalized permutation $\pi$ associated to the linear involution of Figure~\ref{ex:giem} is:
$$
\pi=\left(\begin{array}{ccccc}
1 & 2 & 3 &2&4\\4&5&1&3&5
\end{array}\right). 
$$
\end{exemple}

\begin{remark}
As we have seen before, an interval exchange map can be seen as a linear involution. Also, the table representations of the corresponding combinatorial data are the same. In the sequel, the definitions and statements that we give are valid for linear involutions and for interval exchange maps. 
\end{remark}

\subsection{Rauzy induction and Rauzy classes}
When $T:X\rightarrow X$ is a interval exchange transformation, the first return map of $T$ on a subinterval $X'\subset X$ is still an interval exchange map. The image of $T$ by the \emph{Rauzy induction} $\mathcal{R}$ is the first return map of $T$ on the biggest subinterval $X'\subsetneq X$ which has the same left end as $X$, and such that $\mathcal{R}(T)$ has the same number of intervals as $T$ (see~\cite{Veech82,Marmi:Moussa:Yoccoz}).

Similarly, we can define Rauzy induction for linear involutions by considering first return maps on $X'\times \{0,1\}$, when $X'\subset X$ (see Danthony and Nogueira \cite{DaNo}).

Let $T=(\pi,\lambda)$ be a linear involution on $X$ and denote by $(l,m)$ the type of $\pi$. We identify $X$ with the interval $(0,L)$. If $\lambda_{\pi(l)}\neq \lambda_{\pi(l+m)}$, then the Rauzy induced $\mathcal R(T)$ of $T$ is the linear involution obtained by the first return map of $T$ to
$$
\bigl(0,\max(L-\lambda_{\pi(l)}, L-\lambda_{\pi(l+m)})\bigr) \times
\{0,1\}.
$$
The combinatorial data of the new linear involution depends only on the combinatorial data of $T$ and whether $\lambda_{\pi(l)}>\lambda_{\pi(l+m)}$ or $\lambda_{\pi(l)}<
\lambda_{\pi(l+m)}$. We say that $T$ has type $0$ or type $1$ respectively.
The corresponding combinatorial operations are denoted by $\mathcal{R}_0$ and $\mathcal{R}_1$ correspondingly. Note that if $\pi$ is a given generalized permutation, the subsets $\{T=(\pi,\lambda) ,\ \lambda_{\pi(l)}>\lambda_{\pi(l+m)}\}$ or $\{T=(\pi,\lambda), \lambda_{\pi(l)}< \lambda_{\pi(l+m)}\}$ can be empty because $\pi(l)=\pi(l+m)$ or because the nontrivial relation $\sum_{i=1}^l \lambda_{\pi(i)}=\sum_{j=l+1}^{l+m}\lambda_{\pi(j)}$ that must be fulfilled by $\lambda$.

Let us fix some terminology: given $k\in \{1,\ldots,l+m\}$, the \emph{other occurrence} of the symbol $\pi(k)$ is the unique integer $k'\in \{1,\ldots,l+m\}$, distinct from $k$, such that $\pi(k')=\pi(k)$. In order to describe the combinatorial Rauzy operations $\mathcal R_0$ and $\mathcal R_1$, we first define two intermediary maps $\widetilde{\mathcal R}_0$, $\widetilde{\mathcal{R}}_1$: 

\begin{enumerate}
\item We define $\widetilde{\mathcal R}_0$ in the following way: \\
$\bullet$ If the other occurrence $k$ of the symbol $\pi(l)$ is in $\{l+1,\ldots,l+m-1\}$, then we define $\widetilde{\mathcal R}_0(\pi)$ to be of type $(l,m)$ obtained by removing the symbol $\pi(l+m)$ from the occurrence $l+m$ and putting it at the occurrence $k+1$, between the symbols $\pi(k)$ and $\pi(k+1)$. \\
$\bullet$ If the other occurrence $k$ of the symbol $\pi(l)$ is in $\{1,\ldots,l-1\}$, and if there exists another symbol $\alpha$, whose both occurrences are in $\{l+1,\ldots,l+m-1\}$, then we we define $\widetilde{\mathcal R}_0(\pi)$ to be of type $(l+1,m-1)$ obtained by removing the symbol $\pi(l+m)$ from the occurrence $l+m$ and putting it at the occurrence $k$, between the symbols $\pi(k-1)$ and $\pi(k)$ (if $k=1$,  by convention the symbol $\pi(l+m)$ is put on the left of the first symbol $\pi(1)$).\\
$\bullet$ Otherwise $\widetilde{\mathcal R}_0 \pi$ is not defined. 

\item The map $\widetilde{\mathcal{R}}_1$ is obtained by conjugating $\widetilde{\mathcal{R}}_0$ with the transformation that interchanges the two lines in the table representation.
\end{enumerate}
Then, $\mathcal{R}_0(\pi)$  (resp. $\mathcal{R}_1(\pi)$) is obtained by renumbering $\widetilde{\mathcal{R}}_0(\pi)$ (resp. $\widetilde{\mathcal{R}}_1(\pi)$) to get a reduced generalized permutation. 
For another definition of $\mathcal{R}_0$ and $\mathcal{R}_1$  in terms of the map $\pi$, we refer to \cite{Boissy:Lanneau}.

\begin{Example}
Let us consider the generalized permutation 
$\pi=\left(\begin{smallmatrix} 1 & 2 & 3 & 4 & 3\\ 2 & 4 & 5 & 5 & 1
\end{smallmatrix}\right)$. We have
$$
\widetilde{\mathcal R}_0(\pi)=\left(\begin{array}{cccccc} 1 & 2 & 1 & 3 & 4 & 3\\ 2 & 4 & 5 & 5 &&
\end{array}\right) =
\mathcal R_0(\pi),
$$
 and
 $$
\widetilde{ \mathcal R}_1(\pi)=\left(\begin{array}{cccccc} 1 & 3&2&3&4\\ 2&4&5&5&1
\end{array}\right) \textrm{ so }
\mathcal R_1(\pi)=\left(\begin{array}{cccccc} 1 &2 & 3 & 2 & 4 \\ 3 & 4 & 5 & 5 & 1
 \end{array}\right).
$$
\end{Example}

\begin{figure}[htb]
\includegraphics[scale=0.8]{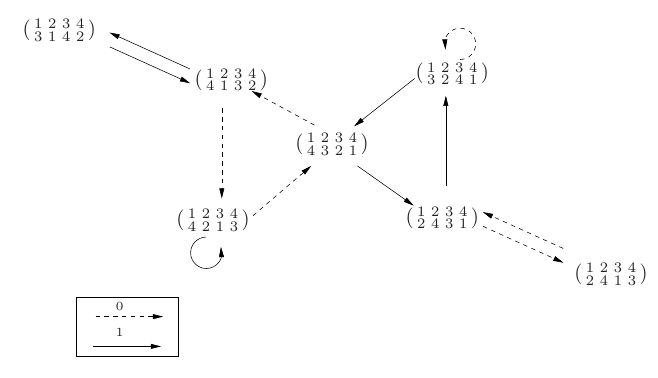}
\caption{An example of a Rauzy diagram for permutations.}
\label{graphe:h11}
\end{figure}

\begin{figure}[htb]
\includegraphics[scale=0.8]{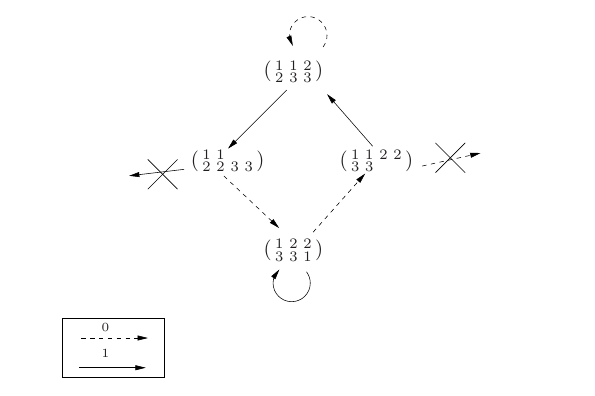}
\caption{An example of a Rauzy diagram for generalized permutations.}
\label{graphe:qiiii}
\end{figure}

\begin{defs}
A \emph{Rauzy class} is a minimal subset of reduced generalized permutations (or permutations) which is invariant by the combinatorial Rauzy maps $\mathcal{R}_0,\mathcal{R}_1$. 
A \emph{Rauzy diagram} is the oriented graph whose vertices are the set of elements of a Rauzy class, and whose edges correspond to the transformations $\mathcal{R}_0$ and $\mathcal{R}_1$.
\end{defs}
\begin{remark}
In this paper, we will speak only of Rauzy class of \emph{irreducible} permutations or generalized permutations (see Definition~\ref{susp:irr} in the next paragraph, and the discussion that follows about irreducibility).
\end{remark}

\subsection{Suspension data and Zippered rectangles construction}
Starting from a linear involution $T$, we want to construct a flat surface $S$ and an  horizontal segment $X$ such that the corresponding first return map of the vertical foliation gives $T$. Such pair $(S,X)$ will be called a \emph{suspension} over $T$, and the parameters encoding this construction will be called \emph{suspension datum}.

\begin{defs} \label{susp:irr}
Let $T=(\pi,\lambda)$ be a linear involution and let
$(\lambda_{k})_{k \in \{1,\ldots,d\}}$ be the lengths of the
corresponding intervals. Let $\{\zeta_{k}\}_{k \in
\{1,\ldots,d\}}$ be a collection of complex numbers such that: 
\begin{enumerate} 
\item $\forall k \in \{1,\ldots,d\} \quad Re(\zeta_{k})=\lambda_{k}$. 
\item $\forall 1\leq i \leq l-1 \quad Im(\sum_{j\leq i} \zeta_{\pi(j)})>0$ 
\item $\forall 1\leq i \leq m-1 \quad Im(\sum_{1\leq j\leq i} \zeta_{\pi(l+j)})<0$
\item $\sum_{1\leq i\leq l} \zeta_{\pi(i)} = \sum_{1\leq j\leq m}\zeta_{\pi(l+j)}$.
\end{enumerate} 
The collection $\zeta=\{\zeta_{i}\}_{i\in \{1,\ldots,d\}}$ is called a \emph{suspension datum} over $T$. The existence of a suspension datum depends only on $\pi$, hence we will say that $\pi$ is \emph{irreducible} if $(\pi,\lambda)$ admits a suspension data.
\end{defs}

 We refer to \cite{Boissy:Lanneau} (Section~3) for a combinatorial criterion of irreducibility for the case when $\pi$ does not correspond to an interval exchange map.

This notion of irreducibility is relevant when we consider Rauzy classes for generalized permutations. Indeed, if $\pi$ is irreducible and if $\pi'$ is in the Rauzy class generated by $\pi$ (\emph{i.e.} the set of descendants of $\pi$ after iterating the combinatorial Rauzy inductions), then $\pi'$ is irreducible and $\pi$ is in the Rauzy class generated by $\pi'$. Therefore, being in the same Rauzy class is then an equivalent relation on the set of irreducible generalized permutations. However, this is not necessarily true if we consider generalized permutations that are not necessarily irreducible: indeed, there exists ''nonirreducible'' generalized permutations whose associated Rauzy class contains irreducible generalized permutations  (see \cite{Boissy:Lanneau}, section~5 and Appendix~A). 

\begin{figure}[htbp]

 \begin{center}
 \includegraphics{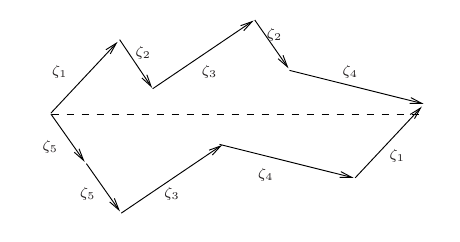}
 \caption{A suspension over a linear involution.}
 \label{figure:suspension:data}
 \end{center}
\end{figure}

Given an interval  exchange map $T$ and a suspension data, there is a well known construction due to Veech, that gives a translation surface and a horizontal segment whose corresponding return map of the vertical geodesic flow is $T$ (see \cite{Veech82,Marmi:Moussa:Yoccoz}). This construction is called the \emph{zippered rectangles construction}.
One can generalize this construction to linear involutions (\cite{B2,Boissy:Lanneau}). Given a suspension datum $\zeta$ over a linear involution $T=(\pi,\lambda)$, we get a flat surface $S$ and an open horizontal segment $X$ (see Figure~\ref{figure:suspension:data}) with an orientation. The first return map of the vertical foliation of $S$ on $X$ is precisely the linear involution $(\pi,Re(\zeta))$. Furthermore, the segment $X$ also satisfies the following properties:
\begin{enumerate}
\item the segment $X$ is adjacent to a singularity on its left,
\item there is a vertical geodesic of $S$ that starts from a singularity and passes through the right end of $X$ before intersecting $X$,
\item any vertical geodesic of $S$ intersects $X$.
\end{enumerate}
We write $(S,X)=Z(\pi,\zeta)$. In fact, the converse is true:

\begin{lem}\label{cs:zipp:rect}
Let $S$ be a flat surface and $X$ be an open horizontal segment $S$ with a choice of orientation. We assume that $X$ satisfies the properties (1)--(2) stated previously, and intersects any vertical saddle connection.
 
 There exists a unique suspension datum $(\pi,\zeta)$, with $\pi$ reduced, such that $(S,X)=Z(\pi,\zeta)$.
\end{lem}

\begin{remark}
In the above lemma, one need $X$ open for technical reasons: it allows us to replace property (3) above, by a condition which is much simpler since there are only a finite number of vertical saddle connections. If $X$ is closed, then $X$ might intersect all vertical saddle connections, but not all vertical geodesics. 
\end{remark}

\begin{proof}
For the case of translation surfaces, the fact that $S$ is obtained by the zippered rectangles construction is a well known fact, and the corresponding permutation and suspension data come from the first return map of the vertical geodesic flow. 
For the case of quadratic differentials, a proof when the surface has no vertical saddle connections can be found in \cite{B2} (Proposition~2.2.). The proof in our case is similar. We give a sketch and refer to \cite{B2} for details.

Let $T=(\pi,\lambda)$ be the linear involution associated to $X$. Up to a finite subset $\Sigma_T$, $X\times\{0,1\}$ is a finite union of open subsets $X_1\ldots,X_{l+m}$, such that $T_{|X_i}$ is a translation or a half-turn. Let $k\neq k'$ be in $\{1,\ldots,l+m\}$ such that $\pi(k)=\pi(k')$. There is an embedded rectangle $R$ whose horizontal edges are identified with $X_k$ and $X_{k'}$. A point in $X$ cannot be in the interior of $R$ since $T$ is the first return map on $X$ of the vertical foliation. Assume that a vertical side of $R$ contains at least two singularities, then it contains a vertical saddle connection, which therefore intersects $X$. Since $X$ is an \emph{open} interval, a subset of $X$ is contained in the interior of $R$, which contradicts the previous assertion.

With this additional argument, one can check that the construction given in \cite{B2}, Proposition~2.2 defines the suspension datum $\zeta$ in a similar way.

\end{proof}

The Rauzy induction on interval exchange maps or on linear involutions admits a natural extension on the space of suspension data. This is called the Rauzy--Veech induction. Let $T=(\pi,\lambda)$ be a linear involution and let $\zeta$ be a suspension data over $T$. We define $\widetilde{\mathcal R}(\pi,\zeta)$ as follows.

\begin{itemize}
\item 
If $T=(\pi,\lambda)$ has type $0$, then
$\widetilde{\mathcal R}(\pi,\zeta)=(\widetilde{\mathcal R}_0\pi,\tilde{\zeta})$, with
$\tilde{\zeta}_k=\zeta_k$ if $k\neq \pi(l)$ and
$\tilde{\zeta}_{\pi(l)}=\zeta_{\pi(l)}-\zeta_{\pi(l+m)}$. 
\item
If $T=(\pi,\lambda)$ has type $1$, then
$\widetilde{\mathcal{R}}(\pi,\zeta)=(\widetilde{\mathcal{R}}\pi,\tilde{\zeta})$, with
$\tilde{\zeta}_k=\zeta_k$ if $k\neq \pi(l+m)$ and
$\tilde{\zeta}_{\pi(l+m)}=\zeta_{\pi(l+m)}-\zeta_{\pi(l)}$. 
\end{itemize} 

Recall that the generalized permutations $\widetilde{\mathcal{R}}_0(\pi)$, $\widetilde{\mathcal{R}}_1(\pi)$ are not necessarily reduced. Hence, after renumerating $\widetilde{\mathcal{R}}(\pi,\zeta)$ in order to get a reduced generalized permuation, we get the pair $\mathcal{R}(\pi,\zeta)$.

\begin{figure}[htbp]
 \begin{center}


 \includegraphics[width=360pt]{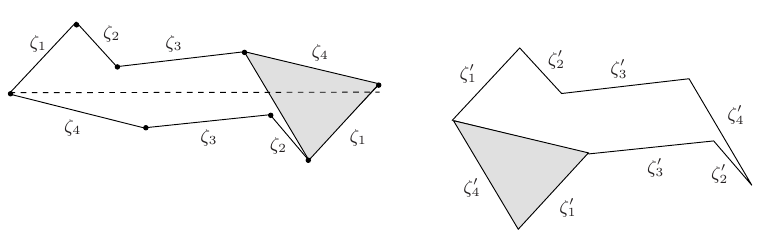}
 \caption{Rauzy-Veech induction on a suspension over an interval exchange transformation.}
 \label{fig:rauzy:suspension}
 \end{center}
\end{figure}

\begin{remark}\label{rauzy:susp}
The pair $\mathcal{R}(\pi,\zeta)=(\pi^{\prime},\zeta^{\prime})$ defines a suspension datum over $\mathcal{R}(T)$. If we denote $(S,X)=Z(\pi,\zeta)$ and $(S', X')=Z(\pi',\zeta')$, the two flat surfaces $S$ and $S'$ are naturally isometric since one can obtain one surface from the other by ``cutting and pasting'' (see Figure~\ref {fig:rauzy:suspension}). Also, under this identification, we have $X'\subset X$.
\end{remark}

Let $\pi$ be a permutation or a generalized permutation and let $\zeta$ be a suspension data. Since the set of suspension data associated to $\pi$ is connected (in fact convex) and the zippered rectangles construction is continuous with respect to the variations of $\zeta$, then all surfaces obtained from a permutation $\pi$ with the zippered rectangles construction belong to the same connected component $\mathcal C(\pi)$ of the stratum.

Let $\mathcal C$ be a connected component of a stratum of the moduli space of Abelian differentials or of quadratic differentials. We denote by $\mathcal{T}_{\mathcal{C}}$ the set 
$$\mathcal{T}_{\mathcal{C}}=\{(\pi,\zeta),\ \mathcal{C}(\pi)=\mathcal{C},\ \zeta \textrm{ is a suspension data for }\pi\},$$ and $\mathcal{H}_{\mathcal{C}}$ the quotient of this set by the Rauzy--Veech induction. The following proposition is clear.

\begin{prop}\label{clear:prop}
The set of connected components of $\mathcal{H}_{\mathcal{C}}$ is in one-to-one correspondence with the set of Rauzy classes $C$ corresponding to a connected component of the moduli space of Abelian or quadratic differentials. 
\end{prop}

\section{Rauzy classes and covering of a stratum}
According to remark \ref{rauzy:susp}, the zippered rectangles construction provides a natural map $\widehat{Z}$ from $\mathcal{H}_\mathcal{C}$ to the ramified covering $\HCC$ of $\mathcal{C}$, obtained by considering the pairs $(S,l)$, where $S\in\mathcal{C}$ and $l$ is a horizontal separatrix adjacent to a singularity of $S$.

\begin{lem}\label{hz:homeo}
The map $\HZ$ is a homeomorphism on its image.
\end{lem}

\begin{proof}
First, let $S$ be such that there exists $(\pi,\zeta)\in \mathcal{T}_{\mathcal{C}}$ with $Z(\pi,\zeta)=(S,*)$. 

We claim that $\zeta$, with the condition $\sum_{i=1}^{l} \zeta_{\pi(i)}=\sum_{j=l+1}^{l+m} \zeta_{\pi(j)}$, defines local coordinates of the ambient stratum. Indeed, in the Abelian case, the sides of the polygon defined by the parameters $\zeta_i$ form a basis of the relative homology, and integrating $\omega$ along this basis gives precisely $\zeta_1,\dots ,\zeta_d$. In the quadratic case, one must consider the natural double cover $\widehat{S}\to S$  , and it is easy to check that integrating the one form corresponding to $\widehat{S}$ along a basis of homology gives $2\zeta$. This proves the claim.  This implies that $Z$ is open, and so is $\HZ$.
 
Now we show that $\HZ$ is injective. 
The pair $(S,l)\in \HCC$ is in the image of $\HZ$ if and only if there exists a segment $X\subset l$, that satisfies the hypothesis of Lemma~\ref{cs:zipp:rect}. For such segment, there  exists a unique $(\pi,\zeta)$ such that $Z(\pi,\zeta)=(S,X)$. 
Now let $X'$ be another such segment, then we must have $X\subset X'$ or $X'\subset X$, and $X'$ defines a new suspension data $(\pi',\zeta')$. We assume for instance that $X'\subset X$. We claim that there exists an integer $n\geq 0$ such that $\mathcal{R}^n(\pi,\zeta)=(\pi',\zeta')$. Assuming the claim, we can conclude that there exists a unique class $[(\pi,\zeta)]\in \mathcal{H}_{\mathcal{C}}$ in the preimage of $(S,l)$ by the map $\HZ$.

When $S$ is a translation surface without vertical saddle connections, the claim is Proposition~9.1 of \cite{Veech82}. We prove the claim in the general case. Let us consider the (possibly finite) sequence of iterates of $(\pi,\zeta)$ by the Rauzy induction. We denote $\mathcal{R}^n(\pi,\zeta)=(\pi^{(n)},\zeta^{(n)})$ 
 and $T^{(n)}$ the corresponding linear involution. We identify the interval $X^{(n)}$ (resp. $X'$) with the interval $]0,x^{(n)}[$ (resp. $]0,x'[$) of $\R$. Three cases are possible. 

\begin{enumerate}
\item There exists $n>0$ such that $x^{(n)}<x'$. We denote by $n_0$ the biggest integer such that $x^{(n_0)}>x'$. By definition of $X'$, there is a vertical geodesic $\gamma$ starting from $x'$ and that hits a singularity before intersecting the interval $]0,x'[$. We claim that it doesn't intersect the interval $]x',x^{(n_0)}[$. Indeed, if $\gamma$ intersects $]x',x^{(n_0)}[$ before hitting a singularity, then we consider $x''\in ]x',x^{(n_0)}[ $ the rightmost intersection point. We must have $x''\leq x^{(n_0+1)}$ which contradicts the hypothesis on $n_0$.

It follows that  $T^{(n_0)}$ is not defined on $(x',\varepsilon)$, for $\varepsilon$ corresponding to the direction of $\gamma$. We know by hypothesis that $\mathcal{R}(\pi^{(n_0)},\zeta^{(n_0)})$ exists, and by definition of the Rauzy induction, we have $x^{(n_0+1)}=x'$. Hence, $(\pi',\zeta')=\mathcal{R}^{(n_0+1)}(\pi,\zeta)$.

\item There exists $n$ such that $x^{(n)}>x'$ and $\mathcal{R}(\pi^{(n)},\zeta^{(n)})$ is not defined. This means that there exists $x^{(n+1)}\geq x'$ such that $T^{(n)}(x^{(n+1)},0)$ and $T^{(n)}(x^{(n+1)},1)$ are not defined. Then there is a saddle connection $\gamma$ that intersects $X^{(n)}$ only in the point $x^{(n+1)}$. Hence, $X'=]0,x'[$ does not intersect $\gamma$, contradicting the hypothesis on $X'$.

\item The sequence $(\pi^{(n)},\zeta^{(n)})$ is infinite and for all $n$, $x^{(n)}>x'$. The sequence $(x^{(n)})_n$ is decreasing and bounded from below.  Hence it converges to a limit $x^{(\infty)}$ which is greater than, or equal to $ x'$. According to the proof of Proposition~4.2 in \cite{Boissy:Lanneau} $T^{(n)}(x^{(\infty)},0)$ and $T^{(n)}(x^{(\infty)},1)$ are not defined for $n$ large enough. Then, there is a saddle connection $\gamma$ that intersect $X^{(n)}$ only in the point $x^{(\infty)}$.  Hence, $X'=]0,x'[$ does not intersect $\gamma$, contradicting the hypothesis on $X'$.
\end{enumerate}

\end{proof}

\begin{prop}
The complement of $\HZ(\mathcal{H}_{\mathcal{C}})$ is contained in a subset of $\HCC$ which is a countable union of real analytic codimension 2 subsets.
\end{prop}

\begin{proof}
 If $S$ has no horizontal saddle connections, any horizontal geodesic is dense.  Hence,  a horizontal segment $X$ adjacent to a singularity will intersect all the vertical saddle connections, as soon as this segment is long enough and by Lemma~\ref{cs:zipp:rect}, the pair $(S,X)$ is in the image of $Z$ for a well chosen $X$. We can also apply Lemma~\ref{cs:zipp:rect} if $S$ has no vertical saddle connection. 
 
Now if $(S,l)\in \HCC$ is such that $S$ has no vertical or no horizontal saddle connections, then $(S,l)$ is in the image of $\HZ$.  Hence,  the complement of the image of $\HZ$ is contained in the set of elements in $\HCC$ whose corresponding flat surface has at least a vertical and a horizontal saddle connections. This set is a countable union of real analytic codimension 2 subsets.
\end{proof}

\begin{cor}\label{hatc:rauzy}
The number of Rauzy classes corresponding to a connected component $\mathcal{C}$ of the moduli space of Abelian or quadratic differentials is equal to the number of connected components of $\HCC$. 
\end{cor}

\begin{proof}
From Proposition~\ref{clear:prop} and Lemma~\ref{hz:homeo}, we just need to prove that the number of connected components of $\HCC$ is equal to the number of connected component of $\HZ(\mathcal{H}_C)$. It is a standard fact that removing a codimension two subset to a smooth manifold does not change its number of connected components. In our case, we remove from an orbifold a countable union of codimension 2 subsets. 

Let $x_1$ and $x_2$ be elements of $\HZ(\HHCC)$ and in the same connected component of $\HCC$. We want to construct a path in $\HZ(\HHCC)$ that joins $x_1$ and $x_2$. Up to considering a local chart of $\HCC$, we can assume that $x_1$ and $x_2$ are in an open subset $\Omega$ of $\C^k$, and there is a finite group $G$ acting on $\Omega$ such that $\Omega/G$ is homeomorphic to an open subset $U$ of $\HCC$. By definition, a real analytic codimension $2$ subset in $U$ corresponds to a real analytic codimension 2 subset of $\Omega$.  Hence,  the elements of $U\backslash \HZ(\HHCC)$ correspond to a countable union $\cup_{i\in \N} F_i$ of smooth codimension 2 subsets of $\Omega$.
Without loss of generality, we can assume that $\Omega$ is convex. Consider a real hyperplane $H$ separating $x_1$ and $x_2$. For each codimension 2 subset $F_i$, the set of elements $y\in H$ such that at least one of the segments $[x_1,y]$ or $[x_2,y]$ contains an element of $F_i$ is of measure zero for the natural Lebesgue measure in $H$. 
 Hence,  the set of elements $y\in H$ such that at least one of the segments $[x_1,y]$ or $[x_2,y]$ intersects $\cup_{i\in \N} F_i$ is of measure zero. So, there is an element $x\in H\cap \Omega$ such that neither $[x_1,x]$ nor $[x,x_2]$ intersects $\cup_{i\in \N} F_i$. This defines a suitable path joining $x_1$ and $x_2$. This concludes the proof.
\end{proof}

\begin{prop}\label{RC:Cm}
The number of distinct Rauzy classes corresponding to a connected component $\mathcal{C}$ of the moduli space of Abelian or quadratic differentials, is equal to the number of connected components of the covering of $\mathcal{C}$ that we obtain by marking a singularity.
\end{prop}

\begin{proof}
Remark that if two separatrices $l_1$ and $l_2$ are adjacent to the same singularity, the two pairs $(S,l_1)$ and $(S,l_2)$ are in the same connected component of $\HCC$, then apply Corollary~\ref{hatc:rauzy}.
\end{proof}

\section{Marked flat surfaces}

In this section, we compute the connected components of the moduli space of flat surfaces with a marked singularity. We will study separately the Abelian and quadratic case.

\subsection{Moduli space of Abelian differentials with a marked singularity}

Here, we assume that $\CC$ is a connected component of the moduli space of Abelian differentials. Recall that the \emph{degree} of a singularity in a translation surface is the integer $k$ such that the corresponding conical angle is $(k+1)2\pi$.

We consider the ramified covering $\mathcal{C}_m$ of $\mathcal{C}$ to be the moduli space of pairs $(S,P)$, where $S\in \mathcal{C}$ and $P$ is a singularity of $S$. According to Proposition~\ref{RC:Cm}, we must count the number of connected components of $\mathcal{C}_m$.

The goal of this section is to prove Proposition~\ref{CC:Cm:ab}, which will complete the proof of Theorem~A for Abelian differentials.

\begin{prop} \label{CC:Cm:ab}
Let $\mathcal{C}$ be a connected component of a stratum in the moduli space of Abelian differentials and let $\mathcal{H}(k_1^{\alpha_1},\ldots,k_r^{\alpha_r})$, with $k_i\neq k_j$ for $i\neq j$, and $k_i\geq 0$ and $\alpha_i>0$ for each $i$, be the ambient stratum. Then $\CC_m$ admits exactly $r$ connected components.
\end{prop}

We want to show that $(S_1,P_1)$ and $(S_2,P_2)$  in $\mathcal{C}_m$ are in the same connected component if and only if the degree of $P_1$ is equal to the degree of $P_2$.  
If $(S_1,P_1)$ and $(S_2,P_2)$ are in the same connected component of $\mathcal{C}_m$, then the degree of $P_1$ is clearly equal to the degree of $P_2$.  
We want to prove the converse. Since $\mathcal{C}_m$ is a ramified covering of $\mathcal{C}$, it is enough to show this for $S_1=S_2$.


%


\medskip

For the following definition, note that a saddle connection persists under any small deformation of the surface inside the ambient stratum. 

\begin{defs}
Let $S$ be a translation surface.
A saddle connection on $S$ is \emph{simple} if, up to a small deformation of $S$ inside the ambient stratum, there are no other saddle connections parallel to it.
\end{defs}

\begin{lem}\label{swap}
Let $S\in \mathcal{C}$ and $P_1, P_2$ be two singularities of the same degree. If there exists a simple saddle connection between $P_1$ and $P_2$, then $(S,P_1)$ and $(S,P_2)$ are in the same connected component of $\mathcal{C}_m$. 
\end{lem}

\begin{proof}
We denote by $\gamma$ the simple saddle connection between $P_1$ and $P_2$, and by $k$ the degree of $P_1$ and $P_2$. We can also assume that $\gamma$ is vertical and up to a slight deformation of $S$, there is no saddle connections parallel to $\gamma$. Recall that the Teichmüller flow acts continuously, so we can apply to $S$ the Teichmüller geodesic flow, and obtain a surface surface $S'=g_tS$ in the same connected component as $S$. There is a natural bijection from the saddle connections of $S$ to the saddle connections of $g_tS$. The holonomy vector $v=(v_1,v_2)$ of a saddle connection becomes $v_t=(e^{t/2} v_1,e^{-t/2} v_2)$. This implies that the length of a given saddle connection in $S'$ divided by the length of $\gamma'$ corresponding to $\gamma$ tends to infinity, as $t$ tends to infinity. 
The set of holonomy vectors of saddle connections is discrete, and therefore, if $t$ is large enough, we can assume that the saddle connection $\gamma'$ is very small compared to any other saddle connection of $S'$. The two singularities corresponding to $P_1$ and $P_2$, that we denote by $P_1'$ and $P_2'$, are the endpoints of $\gamma'$. It is sufficient to show that $(S',P_1')$ and $(S',P_2')$ are in the same connected component of~$\mathcal{C}_m$. If $t$ is large enough, then $S'=g_t.S$ is obtained after breaking up a zero of degree $2k$ into two zeroes of degree $k$, using the local construction described in section~\ref{bzero:loc}.

\medskip

%

The small saddle connection that appear in the procedure corresponds to $\gamma'$. 
 In this procedure, we can continuously turn the parameter defining $\gamma'$, and therefore $(S',P_1')$ and $(S',P_2')$ are in the same connected component of $\mathcal{C}_m$ (see Figure~\ref{tourne:zeros}). 

\begin{figure}[htb]
\begin{center}
\includegraphics[width=360pt]{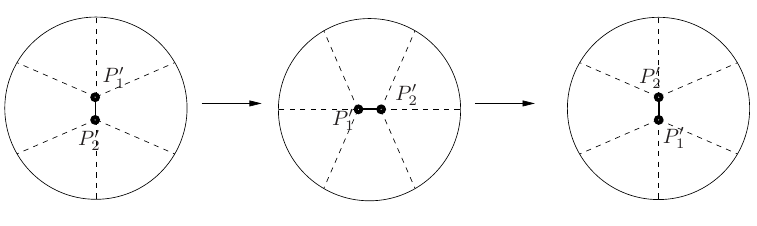}
\caption{Interchanging two zeros of the same degree.}
\label{tourne:zeros}
\end{center}
\end{figure}

\medskip
\end{proof}

Now given a flat surface $S\in \CC$ and two singularities $P,Q$ of the same degree, one would like to find a simple saddle connection that joins $P$ and $Q$. In fact, it is enough to find a broken line that consists of simple saddle connections whose endpoints are singularities of the same degree as $P$ and $Q$. This is the main idea of the proof of Proposition~\ref{CC:Cm:ab}.

\medskip

\begin{proof}[Proof of Proposition~\ref{CC:Cm:ab}]
For each $k$, we show that the subset of $\CC_m$ corresponding to a singularity of degree $k$ is connected. For this, it is enough to find a surface $S\in \mathcal{C}$, and a collection of simple saddle connections connecting all the singularities of degree $k$. Without loss of generality, we assume that $k=k_1$.

We use the following construction: we start from a surface $S_0\in \mathcal{H}(\alpha_1 k_1,k_2^{\alpha_2},\ldots,k_r^{\alpha_r})$. Then, we break up the singularity of degree $\alpha_1 k_1$ into a singularity of degree $k_1$ and a singularity of degree $(\alpha_1-1)k_1$. We get a surface $S_1\in \mathcal{H}(k_1,(\alpha_1-1) k_1,k_2^{\alpha_2},\ldots,k_r^{\alpha_r})$, and a small simple saddle connection between a singularity $P_1$ of degree $k_1$ and a singularity $Q_1$ of degree $(\alpha_1-1)k_1$. Then, we break up the singularity $Q_1$ into a singularity $P_2$ of degree $k_1$ and a singularity $Q_2$ of degree $(\alpha_1-2)k_1$. There is a simple saddle connection between $P_2$ and $Q_2$, if we choose well our breaking procedure, and if the newly created saddle connection is small enough, then the saddle connection between $P_1$ and $P_2$ persists.

Iterating this process, we finally get a surface $S$ in $\mathcal{H}(k_1^{\alpha_1},k_2^{\alpha_2},\ldots,k_r^{\alpha_r})$ and $P_1,\ldots,P_{\alpha_1}$ with a saddle connection $\gamma_i$ between $P_i$ and $P_{i+1}$, for all $1\leq i\leq \alpha_1-1$. Moreover, all the singularities $P_i$ and the corresponding saddle connections $\gamma_i$ are in a flat disk $D$. Each $\gamma_i$ can be assumed to be very short compared to any other saddle connection which is not entirely in $D$. Now assume that one of the saddle connection $\gamma_i$ is not simple. Then, up to a small deformation of $S$, there is another saddle connection $\gamma_i'\subset D$ which is homologous to $\gamma_i$.  Hence,  $\gamma_i$ and $\gamma_i'$ are the boundary of a metric disk $D'\subset D$. The boundary of $D'$ consists of two parallel saddle connections of the same length. Therefore, we can glue them together by a suitable isometry, and obtain a flat sphere that contains at most two poles that correspond to the end points of $\gamma_i$ and $\gamma_i'$. Such flat sphere cannot satisfy the Gauss-Bonnet equality, which contradicts the fact that $\gamma_i$ is not simple. 

 Hence,  we have proven that our construction provides a surface $S$, with a broken line that consists of a union of simple saddle connections joining all the singularities of degree $k$.  We can apply Lemma~\ref{swap} for each pairs $(P_i,P_{i+1})$, and we get that the $\{(S,P_i)\}_{i\in \{1,\ldots,\alpha\} }$ are in the same connected component of the corresponding moduli space of marked translation surfaces. 
It remains to check that $S$ can be taken in any connected component of $\mathcal{H}(k_1^{\alpha_1},\ldots,k_r^{\alpha_r})$.

Without loss of generality, we can assume that there are no singularities of degree zero, since these degree zero singularities just correspond to regular marked point on the surface, and this is deduced from the other case in a trivial way.

If $S_0$ is in $\HH(2g-2)$, and $S$ is in $\HH(g-1,g-1)$, then $S$ is in the hyperelliptic connected component if and only if  the same is true for $S_0$ (see \cite{Kontsevich:Zorich}).

If $S_0$ is not in the hyperelliptic connected component of $\HH$ and if all the singularities of $S$ have even degree, then breaking up a singularity does not change the parity of the spin structure. Indeed, the breaking procedure does not change the metric outside a small disk and the paths that we choose to compute the parity of spin structure can avoid this disk. Hence,  starting from $S_0$ with even or odd spin structure, we get an even or an odd spin structure.

Therefore, in any connected component $\CC$, there is a surface $S$ obtained by the construction. This proves the proposition.
\end{proof}

\subsection{Moduli space of quadratic differentials with a marked singularity}

\begin{NNrem}
Here, we deal with the moduli space of \emph{quadratic} differentials. Therefore, the \emph{order} of a singularity is the integer $k\geq -1$ such that that the corresponding conical angle is $(k+2)\pi$. Recall that $k=0$ corresponds to a regular marked point on the surface
\end{NNrem}

We want to prove Proposition~\ref{CC:Cm:q}, which will complete the proof of Theorem~A. This proposition is a ``quadratic analogous'' of Proposition~\ref{CC:Cm:ab}.

\begin{prop} \label{CC:Cm:q}
Let $\mathcal{C}$ be a connected component of a stratum in the moduli space of quadratic differentials. Let $\QQ (k_1^{\alpha_1},\ldots,k_r^{\alpha_r})$ be the ambient stratum, with $k_i\neq k_j$ for $i\neq j$, and $k_i\geq -1$ and $\alpha_i>0$. Then $\CC_m$ admits exactly $r$ connected components.
\end{prop}

Although the main ideas of the proof are similar, there are some technical difficulties. For instance, the ``quadratic version'' of Lemma~\ref{swap} is still true, but the proof needs some additional tools. Indeed, the ``singularity breaking up procedure'' introduced in section \ref{bzero:loc} does not work when we break up a singularity of even order into two singularities of odd order. So we must use the non local procedure described in section \ref{sec:non:loc}.

The next two lemma are  ``quadratic'' versions of Lemma~\ref{swap}. Lemma~\ref{swap:quadr} is for singularities of non-negative order and Lemma~\ref{lem:pol} is for poles.
\begin{lem}\label{swap:quadr}
Let $\CC$ be a connected component of a stratum in the moduli space of quadratic differentials.
Let $S\in \CC$ and $P_1, P_2$ be two singularities of the same order $k$, with $k\neq -1$. We assume that there exists a simple saddle connection between $P_1$ and $P_2$. Then $(S,P_1)$ and $(S,P_2)$ are in the same connected component of $\CC_m$.
\end{lem}
\begin{proof}
 When $k$ is even, the proof is exactly the same as in Lemma~\ref{swap}. So we assume that $k$ is odd. As in the proof of Lemma~\ref{swap}, we can assume that the simple saddle connection $\gamma$ of the hypothesis is very small compared to any other saddle connection.


 There exists $S_0$, a path $\nu_0\subset S_0$, and $\varepsilon$ such that $(S,P_1)=\Psi(S_0,\nu_0, \varepsilon)$  (see section~\ref{sec:non:loc} for the definition of the map $\Psi$). Fixing $S_0$, we can make $\varepsilon$ arbitrarily small since $\varepsilon\mapsto \Psi(S_0,\nu_0, \varepsilon)$ is continuous.
 
 Then, we consider a homotopy $(\nu^\theta)_{\theta\in [0,(k+1)\pi]}$, such that $\nu^0=\nu_0$, and $\nu^{\theta}$ is a polygonal curve transverse to the foliation $\mathcal{F}_\theta$ in a neighborhood of $P$. The map $\theta\mapsto \Psi_\theta(S_0,\nu^\theta, \varepsilon)$ is well defined and continuous for $\varepsilon$ small enough. This way, we get a surface $\Psi(S_0,\nu_1, \varepsilon)$. The path $\nu_1$ starts from the sector $I\!I$ and ends in the sector $I$ of $P$. It is natural to compare $\nu_1$ with $\nu_0^{-1}$ (\emph{i.e.} $\nu_0$ with reverse orientation), but these two paths are \emph{a priori} very different (see Figure~\ref{tourne:non:loc}). 

\begin{figure}[htb]
\begin{center}
\includegraphics[width=360pt]{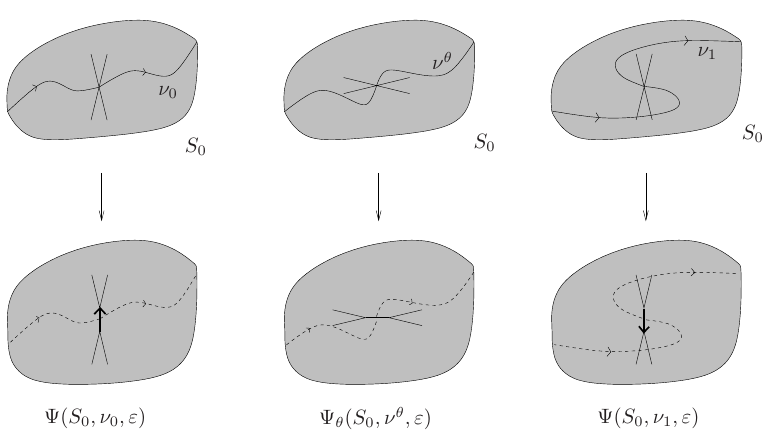}
\caption{Interchanging two singularities of odd order}
\label{tourne:non:loc}
\end{center}
\end{figure}

Using the results stated in section~\ref{sec:non:loc}, there exists $S_1$ in a neighborhood of $S_0$ such that $\Psi(S_0,\nu_1,\varepsilon)=\Psi(S_1,\nu_0^{-1},\varepsilon)$. The surface $S_1$ can be arbitrarily close to $S_0$ as soon as $\varepsilon$ is small enough. Then, we choose a small path joining $S_1$ and $S_0$, and we get therefore a path joining $\Psi(S_1,\nu_0^{-1}, \varepsilon)$ to $\Psi(S_0,\nu_0^{-1}, \varepsilon)$. 

 Hence,  we have built a path joining $\Psi(S_0,\nu_0, \varepsilon)$ to $\Psi(S_0,\nu_0^{-1}, \varepsilon)$. The first (marked) surface is $(S,P_1)$ while the second one is $(S,P_2)$. The lemma is proven. 
\end{proof}

A surface in $\CC_m$ might contain poles. The previous lemma does not work if the marked point is a pole. We need the following:

\begin{lem}\label{lem:pol}
Let $\CC$ be a connected component of a stratum in the moduli space of quadratic differentials.
Let $S\in \CC$ and $P_1, P_2$ two poles. We assume that there exists a saddle connection between $P_1$ and $P_2$. Then $(S,P_1)$ and $(S,P_2)$ are in the same connected component of $\CC_m$.
\end{lem}
\begin{proof}
The saddle connection $\gamma$ joining $P_1$ and $P_2$ is never simple. Indeed, $P_1$ and $P_2$ are in the boundary of a cylinder whose waist curves are parallel to $\gamma$. One side of this cylinder consists of $\gamma$, the opposite side is a union of saddle connections that are necessary parallel to $\gamma$. So $\gamma$ cannot be simple.

In this case, $(S,P_1)$ and $(S,P_2)$ can be joined by performing a suitable Dehn twist on the corresponding cylinder.
\end{proof}

Now we have the necessary tools to prove Proposition~\ref{CC:Cm:q}.

\begin{proof}[Proof of Proposition~\ref{CC:Cm:q}]
We must show that the subset of $\CC_m$ that corresponds to surfaces with a marked point of order $k$, where $k$ is a fixed element of $k_1,\ldots,k_r$ is connected. Without loss of generality, we can assume that $k=k_1$. Also, we can assume that all $k_i$ are nonzero. 

First we assume that $k_1=-1$. 
According to Lanneau (\cite{Lanneau:cc}), there is a surface $S$ in $\CC$ whose horizontal foliation consists of one cylinder. This means we can present   such surface as a rectangle with the following indentifications on its boudary:  
\begin{itemize}
\item the two vertical sides are identified by a translation,
\item the horizontal sides admit a partition of segments which come by pairs of segments of the same length
\item for each such pair, we identify the corresponding segments by translation or by a half-turn. 
\end{itemize}

We can also assume that the corners of the rectangle correspond to singularities. Now, let $P_1$ and $P_2$ be two singularities of order $-1$. Each pole corresponds to two adjacent segments that are identified with each other by a half-turn. 
If these two singularities are on opposite sides of the rectangle, then we get a saddle connection joining $P_1$ and $P_2$ by considering the line joining $P_1$ and $P_2$ in the rectangle. If $P_1$ and $P_2$ are in the same side of the rectangle, then we can slightly deform the corresponding segments in the 1-cylinder decomposition, and this way join the two poles $P_1$ and $P_2$ by a saddle connection (see Figure~\ref{poles}). In any case, we have the desired result (when $k=-1$) in view of Lemma \ref{lem:pol}.

\begin{figure}[htb]
\includegraphics[width=250pt]{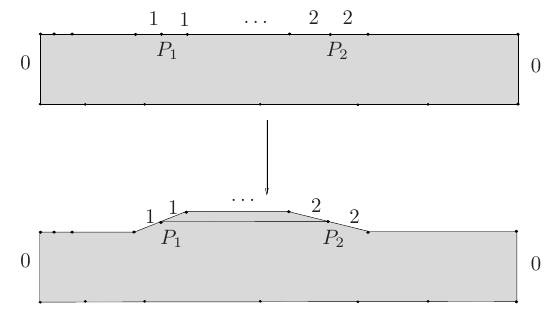}
\caption{Interchanging two poles on a surface with a~1-cylinder decomposition.}
\label{poles}
\end{figure}

Now we assume that $k_1\neq -1$. We first explain the general construction.
By a similar argument as in Proposition~\ref{CC:Cm:ab}, we start from a surface $S_0$ with a singularity $P$ of order $\alpha_1 k_1$ and we break up this singularity into $\alpha_1$ singularities $P_1,\ldots,P_{\alpha_1}$ of order $k_1$. There is a collection of saddle connections joining $P_i$ to $P_{i+1}$ for each $1\leq i\leq \alpha_1-1$. We can assume that $P_1,\ldots,P_{\alpha_1}$ are in a small metric disk $D$. 
Now assume that one of the saddle connection $\gamma_i$ is not simple. Then, up to a slight deformation of $S$, there is another saddle connection $\gamma_i'$ parallel to $\gamma_i$, such that $S\backslash \bigl( \gamma_i \cup \gamma'_i\bigr)$ admits a connected component with trivial linear holonomy (since $\gamma_i$ and $\gamma'_i$ are \^homo\-logous, see~\cite{Masur:Zorich}, Proposition~1 and Theorem~1). However, since $S_0$ has nontrivial linear holonomy, $S\backslash D$ has nontrivial linear holonomy too. Hence, $\gamma_i$ and $\gamma'_i$ are the boundary of a small metric disk $D'\subset D$, which is a contradiction. However, as we will see, we cannot reach any connected component $\CC$ in this way.
\medskip

1- We first assume that the stratum $\QQ=\QQ(k_1^{\alpha_1},\ldots,k_r^{\alpha_r})$ does not contain a hyperelliptic connected component and is not one of the exceptional stratum. Then our connected component $\CC$ is the whole stratum. 
If we start from an initial flat surface $S_0\in \QQ(\alpha_1 k_1,k_2^{\alpha_2},\ldots,k_r^{\alpha_r})$ and perform the previous construction, we get a surface $S\in \mathcal{C}$ and simple saddle connections joining all its singularities of order $k_1$.
We must check that the stratum $\QQ(\alpha_1 k_1,k_2^{\alpha_2},\ldots,k_r^{\alpha_r})$ is not empty. The only strata that are empty are $\QQ(\emptyset), \QQ(1,-1), \QQ(3,1)$ and $\QQ(4)$. Hence, we must have $\QQ\neq \QQ(2,2)$ and $\QQ\neq \QQ(1,1,1,1)$. But these two strata consist only of hyperelliptic flat surfaces, hence $\QQ$ is not one of them by assumption. 
Therefore, using Lemma~\ref{swap:quadr}, we see that $(S,P)$ is in the same connected component of $\CC_m$ for any singularity $P$ of order $k_1$.
\medskip

2- Now we assume that the stratum $\QQ$ is $\QQ(k_1,k_1,k_2,k_2)$, with $k_1\neq k_2$, or $\QQ(k_1,k_1,2k_2+2)$. This stratum has one or two connected components, one of them being hyperelliptic. One can show that in each connected component, on almost every surface $S$, there are simple saddle connections joining the singularities of order $k_1$. (see \cite{B1}, Theorem~3.1 in the case of the hyperelliptic component and \cite{B1} Lemma~4.1 for the other component), and by Lemma~\ref{swap:quadr} we are done. If $\QQ=\QQ(2k_1+2,2k_2+2)$ with $k_1\neq k_2$, there is nothing to prove.
\medskip

3- Assume that $\QQ=\QQ(-1,3,3,3)$. This stratum has two connected components $\CC_{red}$ and $\CC_{irr}$. 
If we start from $S_0\in \QQ(-1,9)$ and break up the singularity of order 9 into three singularities of order 3 as explained previously, we obtain either a surface in $\CC_{red}$ or a surface in $\CC_{irr}$ depending in which connected component we start (see Lanneau~\cite{Lanneau:cc}) and conclude as previously. If the stratum $\QQ$ is one of the other exceptional strata, there is nothing to prove.
\medskip
 
4- We assume that $\QQ=\QQ(k,k,k,k)$.
Let $\mathcal{C}$ be the hyperelliptic connected component of $\QQ$ and $S\in \mathcal{C}$. We denote by $P_{1,1}$,$P_{1,2}$, $P_{2,1}$ and $P_{2,2}$ the singularities of $S$, such that the hyperelliptic involution $\tau$ interchange $P_{i,1}$ and $P_{i,2}$ for $i\in \{1,2\}$. Suppose there is a saddle connection $\gamma$ joining $P_{1,i}$ to $P_{2,j}$ for some $i,j$. Then, $\tau(\gamma)$ is distinct from $\gamma$ and is parallel to $\gamma$, even after a small deformation of $S$. Therefore $\gamma$ is not simple. Hence, $S$ is not obtained from $\QQ(4k)$ by breaking up the singularity as before. 

We can assume that $k\neq -1$, since the other case was already studied. There is a one-to-one mapping from $\CC$ to $\QQ(k,k,-1^{2k+4})$.  Hence, $\CC_m$ is a covering of $\QQ(k,k,-1^{2k+4})$. There exists a surface $S_0\in \QQ(k,k,-1^{2k+4})$ with a simple saddle connection joining its two singularities $P_1$ and $P_2$ of order $k$. We can assume that $S$ is the double covering of $S_0$ ramified over the poles, and that the singularities corresponding to $P_i$ are  $P_{i,1}$ and $P_{i,2}$. For each $i$, there is a simple saddle connection joining $P_{i,1}$ and $P_{i,2}$ (see case (2)), hence the two marked surfaces $(S,P_{i,1})$ and $(S,P_{i,2})$ are in the same connected component of $\CC_m$.
Now we start from $(S,P_{1,1})\in \CC_m$. The corresponding marked surface in $\QQ(k,k,-1^{2k+4})$ is $(S_0,P_1)$. We then consider a path joining $(S_0,P_1)$ and $(S_0,P_2)$ and can lift it to a path joining $(S,P_{1,1})$ to $(S,P_{2,k})$, for some $k\in \{1,2\}$.  Hence, $(S,P_{1,1})$ and $(S,P_{2,1})$ are in the same connected component of $\CC_m$. This proves that $\CC_m$ is connected. 

Let $\mathcal{C}$ be the nonhyperelliptic connected component of $\QQ(k,k,k,k)$. The classification of connected components by Lanneau implies that $k\geq 2$. Then, starting from $S_0\in \QQ(4k)$ and breaking up the singularity into four singularities of degree $k$ as before gives a surface $S\in \mathcal{C}$, since it cannot be in the hyperelliptic connected component as explained before. Hence $\mathcal{C}_m$ is connected.
\medskip

5- If $\QQ=\QQ(2k+2,2k+2)$, the proof is analogous to the previous case.
\end{proof}

\appendix

\section{Computation of the connected component associated to a permutation}

Corollary $B$ states that two irreducible permutations are in the same Rauzy class if and only if the degree of the singularity attached on the left in the Veech construction is the same, and if they correspond to the same connected component of the moduli space of Abelian differentials. 

The first invariant is very easy to compute combinatorially. We give here references for the second invariant.
\begin{itemize}
\item The parity of the spin structure can be computed explicitly from the permutation. This is explained in the paper of Zorich \cite{Zorich:jmd}, Appendix~C. One can also find in Zorich's webpage\footnote{http://perso.univ-rennes1.fr/anton.zorich/} some \emph{Mathematica} program that compute explicitely this invariant.
\item It is strangely not obvious to see whether a permutation corresponds to a hyperelliptic connected component or not. However, in each Rauzy class, we can find a permutation $\pi$ such that $\pi(1)=d$ and $\pi(d)=1$, where $d$ is the number of intervals of the corresponding interval exchange. Such permutation is called cylindrical since it appears naturally for flat surfaces with a one-cylinder decomposition. This was first proven by Rauzy \cite{Rauzy}, but we can find a more constructive proof in \cite{Kontsevich:Zorich}, Appendix~A.3. It is easy to see that the associated connected component is the hyperelliptic one if and only if $\pi$ is the permutation $\pi(k)=d+1-k$. 
Such permutation $\pi$ can be build from another permutation after at most $d^2$ steps of the Rauzy induction in an explicit way (see \cite{Kontsevich:Zorich}).
\end{itemize}

For the case of quadratic differentials, the nonconnected strata are the ones that contain hyperelliptic connected components and the exceptionnal ones. In this case, there is no simple way to decide if two generalized permutations are in the same Rauzy class. 
\begin{itemize}
\item An analogous of the cylindrical permutations exists in each Rauzy classes of generalized permutations, but there is no explicit combinatorial way to find it starting from a given generalized permutation.
\item For the four exceptionnal strata, the only known proof of their nonconnectedness is the explicit computation of the corresponding (extended) Rauzy classes.
\end{itemize}

For related work,  see the paper of Fickenscher \cite{Fick}.

\end{document}